\documentclass{article}
\usepackage{racc}

\usepackage[    
protrusion=true, 
letterspace=50,  
shrink=40,      
factor=1000]    
{microtype}
\usepackage[margin=2.5cm,tmargin=2cm]{geometry}

\let\emphOld\emph

\title{Bijection: Parking-like structures and Tree-like structures}
\author{Jean-Baptiste \textsc{Priez}\thanks{
    \href{www.lri.fr}{Laboratoire de Recherche en
    Informatique}\newline\indent\;\;  
    Bât 650 Ada Lovelace, Université Paris-Sud\newline\indent\;\;
    91405 Orsay, France\newline\indent\;\;
    \href{mailto:jean-baptiste.priez@lri.fr}{jean-baptiste.priez@lri.fr}
    / \href{http://kerios.fr}{kerios.fr}}}

\hypersetup{
    colorlinks = true,
    linkcolor = red!80!black,
    citecolor = gray,
    urlcolor = green!60!black,
    anchorcolor = pink
}
\definecolor{orange}{RGB}{232,101,91}

\definecolor{background}{RGB}{255, 255, 255}
\definecolor{fontcolor}{rgb}{0,0,0}
\everymath{\color{blue!80!black}}
\renewcommand{\emph}[1]{\textcolor{red!70!black}{\emphOld{#1}}}

\let\textOld\text
\renewcommand{\text}[1]{\textOld{\color{fontcolor}#1}}

\pagecolor{background}
\color{fontcolor}
\begin{document}
\color{fontcolor}
\maketitle
\begin{tikzpicture}[remember picture,overlay]
\node [xshift=11cm,yshift=-6.8cm] at (current page.north west)
        [text width=17cm]{
    \small{\textsc{(enumeration, bijection, parking functions, trees,
    species, labeled structures, Lagrange inversion)}} };
\end{tikzpicture}
\hrule\ \\[10pt]
%
    \indent
The \emph{parking problem} is introduced in \cite[\S6]{konheim1966occupancy}:
\begin{quote}
    A car occupied by a man and his dozing wife enters the street on the left
    and moves towards the right. The wife awakens at a capricious moment and
    orders her husband to park immediately! He dutifully parks at his present
    location, if it is empty, and if not, continues to the right and parks at
    the next available space.
\end{quote} 
This \textit{occupancy problem} introduce the \emph{parking function}
combinatorial structures: Let $U$ be a finite set. A parking function $\pf :U
\to \NN_{>0}$ is a function such that,
\begin{align}
     \card \left\{\,\pf(i) \ieq k \;|\; i \in U\,\right\} \seq k\,,&& \text{for
     any } k \ieq u\,.
     \label{eq:pf_condition}
\end{align}
As a hash problem a ``car'' is a key $u \in U$, the ``dozing wives'' act as the
\emph{hash map} $\pf$: they compute a ``park location'' (an index) into the
``street'' (an array of buckets) $\pf(u)$. The ``husbands'' act as an
effective \emph{open addressing hash table} which solve the \emph{hash
collisions} if several wives have same ``caprice''.
In the sequel we keep the wives activities/hash map and we replace the
husbands/open addressing implementation with some other ``implementation'' used
to solve hash collisions. The goal is to enumerate all configurations
associated to a fixed implementation; for example how many \emph{chained hash
tables with linked lists} it is possible to obtain using \emph{parking
functions} as a hash map?\medskip

In \cite{priezvirmaux} the authors redefine (generalized) parking functions
species recursively: Let $\chi: \NN_{>0} \to \NN$ be an non-decreasing function,
\begin{align}
    \PF_\chi &:= \left(\setSp^{\chi(1)}\right)_{0} + 
                 \sum_{n \seq 1} \left(\setSp^{\chi(1)}\right)_{n} \cdot
                                 \PF_{\rho_n}\,,
&&\text{with}\qquad
    \rho_n : m \longmapsto \chi(n + m) - \chi(1)\,,
    \label{eq:pf_def_rec} 
\end{align}
with $\setSp$ the set species; $\setSp^{k}$ denotes the exponentiation of set
species: $\setSp^k := \setSp \cdot \setSp^{k-1}$ with $\cdot$ the
product operator of species. In other terms, $\setSp^k$ is the species of
$k$-sequences of sets ($\setSp^k[U]$ is the sequences of length $k$ of
disjoint sets labeled by $U$).
The notation $+$ is the sum of species and $(\setSp^k)_n$ is the restriction of
the species $\setSp^k$ to sets of cardinality $n$ (the reader may refer to
\cite{bergeron1998combinatorial} about Species Theory or similary 
\cite[chapter II]{flajolet2009analytic} about Labeled Structures).
The parking functions $\PF$ (\ref{eq:pf_condition}) are generalized parking
functions $\PF_\chi$, where we take $\chi$ to be the identity map $Id$ (see
\cite{stanley2002polytope, kung2003goncarov} about generalized parking functions).
\medskip

Despite the fact that the recursive definition (\ref{eq:pf_def_rec}) gives an
efficient way to generate (generalized) parking functions, the underlying
exponential generating series is an inefficient enumerating formula, since the
coefficients $\gf(n) = (n+1)^{n-1}$ are defined as a sum over compositions of
$n$ \cite[Theorem 3.4]{priezvirmaux}.
An elegant way to enumerate \emph{parking functions} is to use bijection with
\emph{forest of rooted trees} $\F$. The species of forest rooted trees $\F$ is
defined by the following recursive definition:
\begin{align}
    \F := \setSp(\singSp \cdot \F)\,,
    \label{eq:species_forest_rooted_trees}
\end{align}
with $\singSp$ the singleton species ($\singSp[U] = \{U\}$ \textit{iff} $U =
\emptyset$) and $F(G)$ denoting the (partitional) \textit{composite} of the
species $G$ in the species $F$.
This bijection of structures gives the generating series of parking functions:
\begin{align*}
    \PF(t) = \sum_{n \seq 0} \gf(n) \frac{t^n}{n!}\,,\qquad
\text{which satisfies the equation:}\qquad
    \PF(t) = \exp(t\,\PF(t))\,.
\end{align*} 
We further generalize the parking functions by replacing
the set species $\setSp$, in (\ref{eq:pf_def_rec}), with other species $\SP$
whose structures model implementations of the hash table (using
\emph{parking functions} as hash maps): the \emph{parking-like structures}. In
the same way we replace $\setSp$ by $\SP$ in the forest definition
(\ref{eq:species_forest_rooted_trees}) and consider the \emph{tree-like
structures}.
We show and make explicit a bijection between the \textit{parking-like
structures} and the \textit{tree-like structures} (Theorem \ref{thm:iso}). From
this isomorphism we obtain generating series defined by functional equation
(Corollary \ref{cor:enum_pf_like}) computable by \emph{Lagrange
inversion}.\medskip

The generic bijection gives interesting correspondance between
\emph{parking-like structures} and trees-like structures: (labeled) \emph{binary
trees}, \emph{$k$-ary trees}, \emph{hypertrees}, \etc. Particularly
the labeled binary trees are isomorphic to the structures of possible
\emph{chained hash tables with linked lists}. 
\section{Parking like species}
    Let $\SP$ be a species.  
Without loss of generality we can suppose it is possible to fully recover the
underlying set $U$ of any $\SP$-structures on $U$. Indeed the species $\bar\SP
:= \SP \times \setSp$ (with $\times$ the cartesian product of species
or equivalently $\bar\SP[U] := \{(\sp, U)$, for any $ \sp \in \SP[U] \}$) is
isomorphic to $\SP$.
%
We also suppose a total order on $U$, noted $<_U$:
%
\begin{align*}
    u_1 <_{U} u_2 <_U \cdots <_U u_n\,,&& \text{with}\; n := \card{U}\,.
\end{align*}
\subsection{Generalization of the generalization}
The definition (\ref{eq:pf_def_rec}) is a constructive definition of the
generalized parking functions \cite{priezvirmaux, stanley2002polytope,
kung2003goncarov}. This equation resumes parking functions on $U$ as a sequence
of sets $(Q_i)$ which satisfies:
\begin{align}
    \sum_{i=1}^{k} \card{Q_i} \seq k\,, &&\text{for any}\; k\ieq n :=
    \card{U}\,.
    \label{eq:pf_condition_sets}
\end{align}
This condition (\ref{eq:pf_condition_sets}) is a translation of the parking
condition (\ref{eq:pf_condition}).  We generalize this species by replacing
$\setSp$ by other species $\SP$:
\begin{align}
    \shift{\SP}{\chi} &:= \left(\SP^{\chi(1)}\right)_{0} + 
                 \sum_{n \seq 1} \left(\SP^{\chi(1)}\right)_{n} \cdot
                                 \shift{\SP}{\rho_n}\,,
    \label{eq:def_shift}
\end{align}
with $\rho_n$ defined as in (\ref{eq:pf_def_rec}).\\

\noindent
\textbf{Notation :} 
We denote $\SP^\btr := \shift{\SP}{Id}$ the parking-like species over $\SP$
associated to the identity map (with $\chi = Id$).\\
 
We call $\shift{\SP}{\chi}$ the \emph{parking like species} over the species
$\SP$ associated to the non-decreasing map $\chi$. (So $\PF_\chi =
\shift{\setSp}{\chi}$ and $\PF = \setSp^\btr$.) In the sequel we are
focusing on this species and we are looking to enumerated
$\SP^\btr$-structures \textit{via} a bijection with \emph{tree-like
structures}.
\begin{prop}
    \label{prop:pf_g_pf_set}
    Let $U$ be a finite set,
    \begin{align*}
    \shift{\SP}{\chi}[U] = \left\{
        (\sp_i)_{i \in [\chi(u+1)]} \;\mid\; \sp_i \in \SP[V_i] \;\;\text{such
        that}\;\; (V_i) \in \PF_\chi[U] 
    \right\}\,.
    \end{align*}
\end{prop}
\begin{remarq}
    The $\shift{\SP}{\chi}$-structures on $U$ are $\chi(u+1)$-sequences of 
    $\SP$-structures such that the $\SP$-structures appearing at indices between
    $\chi(u)+1$ and $\chi(u+1)$ are defined on $\emptyset$.
    If $\chi = Id$ that means the last structure of the sequence is
    a structure on the emptyset. This point will be important in the bijection
    Theorem \ref{thm:iso}.
\end{remarq}
\begin{ex}
    \label{ex:lin_order_pf}
    Let $\linSp$ be the linear order species (defined by $\linSp = \oneSp +
    \singSp \cdot \linSp$ with $\oneSp$ neutral species (such that
    $\oneSp[U] = \{U\}$ \textit{iff} $U = \emptyset$)). The $\linSp$-structures
    are isomorphic to the permutations.
    The $\linSp^\btr$-structures on $[0], [1], [2]$ and $[3]$ are:
    \begin{align*}
    \linSp^\btr[0] &= \left\{(\cdot)\right\}\,; \quad
    \linSp^\btr[1] = \left\{(1 \mid \cdot)\right\}\,; \quad
    \linSp^\btr[2] = \left\{(12\mid\cdot \mid \cdot),\; (21\mid\cdot \mid
    \cdot),\; (1\mid2 \mid \cdot),\; (2\mid1 \mid \cdot)
    \right\}\,;\\
    \linSp^\btr[3] &= \left\{\begin{array}{c}
        (123\mid \cdot\mid \cdot \mid \cdot),\;
        (132\mid \cdot\mid \cdot \mid \cdot),\;
        (213\mid \cdot\mid \cdot \mid \cdot),\;
        (231\mid \cdot\mid \cdot \mid \cdot),\;
        (312\mid \cdot\mid \cdot \mid \cdot),\;
        (321\mid \cdot\mid \cdot \mid \cdot),\\
        (12\mid 3\mid \cdot \mid \cdot),\;
        (21\mid 3\mid \cdot \mid \cdot),\;
        (12\mid \cdot \mid 3 \mid \cdot),\;
        (21\mid \cdot \mid 3 \mid \cdot),\;
        (13\mid 2\mid \cdot \mid \cdot),\;
        (31\mid 2\mid \cdot \mid \cdot),\\
        (13\mid \cdot \mid 2 \mid \cdot),\;
        (31\mid \cdot \mid 2 \mid \cdot),\;
        (23\mid 1\mid \cdot \mid \cdot),\;
        (32\mid 1\mid \cdot \mid \cdot),\;
        (23\mid \cdot \mid 1 \mid \cdot),\;
        (32\mid \cdot \mid 1 \mid \cdot),\\
        (1\mid 23 \mid \cdot \mid \cdot),\;
        (1\mid 32 \mid \cdot \mid \cdot),\;
        (2\mid 13 \mid \cdot \mid \cdot),\;
        (2\mid 31 \mid \cdot \mid \cdot),\;
        (3\mid 12 \mid \cdot \mid \cdot),\;
        (3\mid 21 \mid \cdot \mid \cdot),\\
        (1\mid 2 \mid 3 \mid \cdot),\;
        (2\mid 1 \mid 3 \mid \cdot),\;
        (1\mid 3 \mid 2 \mid \cdot),\;
        (3\mid 1 \mid 2 \mid \cdot),\;
        (2\mid 3 \mid 1 \mid \cdot),\;
        (3\mid 2 \mid 1 \mid \cdot)
    \end{array}
    \right\}\,.
    \end{align*}
\end{ex}
\subsection{Tree-like species and bijection}
Similary to the replacement on the generalized parking functions species, we
replace the set species with another species in the definition of forest species
(\ref{eq:species_forest_rooted_trees}):
\begin{align*}
    \T_\SP := \SP(\singSp \cdot \T_\SP)\;.
\end{align*}
We call $\T_\SP$ the \emph{tree-like species} over the species $\SP$. Focusing
on the identity map, we obtain the following isomorphism with $\SP^\btr$: 
\begin{thm}
    \label{thm:iso}
    There is a bijection between $\SP^\btr$-structures and
    $\T_\SP$-structures ($\SP^\btr[U] \simeq \T_\SP[U]$, for any finite set
    $U$).
\end{thm}
\begin{preuve}[by drawing]
    Let $U$ be a finite set of cardinality $n$.
    Let $\sp$ be a $\SP$-structure on $U$ represented (as in
    \cite{bergeron1998combinatorial}) by:
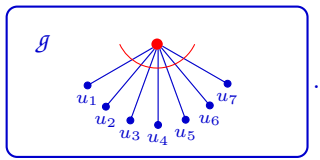
    \begin{align*}
\vcenter{\hbox{\begin{tikzpicture}
\draw[thick, rounded corners] (0,-1) rectangle (4,1);
\fill[red] (2,.5) circle (.5ex);
\foreach \x in {1,2,...,7} {
    \draw (2,.5) ++(-\x*20:.5ex) --++(-\x*20-10:1)[fill] circle (.3ex)
    node[below] {$\scriptstyle u_{\pgfmathparse{int(8-\x)}\pgfmathresult}$};
}
\draw[red] (2.5,.5) arc (-25:-155:.55);
\node at (.5, .5) {$\sp$};
\end{tikzpicture}}}\,.
    \end{align*}

Let $(\sp_i)$ be a $\SP^\blacktriangleright$-structure on $U$ represented by:
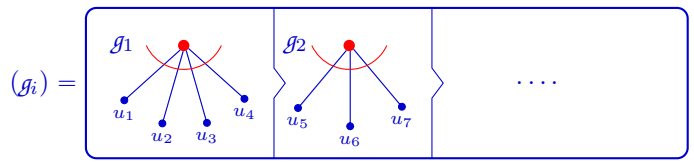
\begin{align*}
    (\sp_i) =
    \vcenter{\hbox{\begin{tikzpicture}
\draw[thick, rounded corners] (0,-1) rectangle (8,1);
\fill[red] (1.3,.5) circle (.5ex);
\foreach \x in {1,2,...,4} {
    \draw (1.3,.5) ++(-\x*32:.5ex) --++(-\x*32-10:1)[fill] circle (.3ex)
    node[below] {$\scriptstyle u_{\pgfmathparse{int(5-\x)}\pgfmathresult}$};
}
\draw[red] (1.8,.5) arc (-25:-155:.55);
\node at (.5, .5) {$\sp_1$};
    \draw (2.5,-1) -- (2.5,-.2) -- (2.65, 0) -- (2.5,.2) -- (2.5,1);
\node at (2.8, .5) {$\sp_2$};
\fill[red] (3.5,.5) circle (.5ex);
 \foreach \x in {1,2,3} {
     \draw (3.5,.5) ++(-\x*40:.5ex) --++(-\x*40-10:1)[fill] circle (.3ex)
     node[below] {$\scriptstyle u_{\pgfmathparse{int(8-\x)}\pgfmathresult}$};
 }
 \draw[red] (4,.5) arc (-25:-155:.55);
    \draw (4.6,-1) -- (4.6,-.2) -- (4.75, 0) -- (4.6,.2) -- (4.6,1);
    \node at (6, 0) {.\,.\,.\,.};
\end{tikzpicture}}}\,,
\end{align*}
This structure is a $n+1$-sequence of $\SP$-structures $\sp_i$ which
satisfies:
\begin{align*}
    \sum_{i=1}^{k} \card{V_i} \seq k&& \text{for any }\; k\ieq n\,;
    \tag{Proposition \ref{prop:pf_g_pf_set}}
\end{align*}
with $V_i$ the underlying set of $\sp_i$, for any $i$ (analogous to
(\ref{eq:pf_condition})).
Using the order $<_U$ on $U$ we associate to the parking like structure
$(\sp_i)$ the total order $<_q$ defined by:
\begin{align*}
      u <_q u' \quad \Longleftrightarrow \quad
      \begin{dcases*}
        u <_U u' & \text{with $u,u'$ elements of the underlying set of
        $\sp_i$,}\\
        j < k & \text{with $u$ (\resp $u'$) an element of the underlying set of
        $\sp_j$ (\resp $\sp_k$),}
      \end{dcases*}
\end{align*}
Following the idea of the construction between parking function and forest of
rooted trees we associate to the $\SP^\blacktriangleright$-structure $(\sp_i)$ a
$\T_\SP$-structure $f$ defined by:
\begin{itemize}
  \item set $\sp_1$ be the root, 
  \item set vertices $u_i \to \sp_{i+1}$, for any $u_i \in U$.
\end{itemize}
This construction is summarized by the following schema (where vertices are
represented by circles):
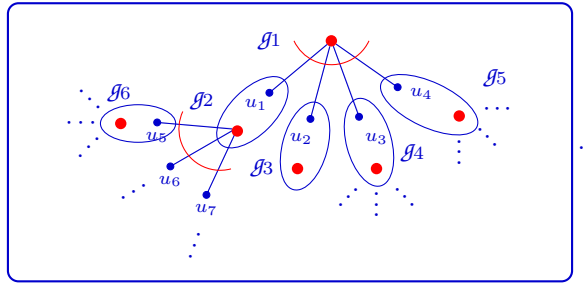
\begin{align*}
 \vcenter{\hbox{\begin{tikzpicture}
\draw[thick, rounded corners] (-3,1) rectangle (4.5,-2.7);
\fill[red] (1.3,.5) circle (.5ex);
\foreach \x in {1,2,...,4} {
    \draw (1.3,.5) ++(-\x*35:.5ex) --++(-\x*35:1)[fill] circle (.3ex); 
}
\draw[red] (1.8,.5) arc (-25:-155:.55); 
\node at (.5, .5) {$\sp_1$};
\node (u1) at (.3, -.35) {$\scriptstyle u_1$}; 
\node (u2) at (.9, -.8) {$\scriptstyle u_2$}; 
\node (u3) at (1.9, -.8) {$\scriptstyle u_3$};
\node (u4) at (2.5, -.2) {$\scriptstyle u_4$};  
\node (g2) at (-.4, -.3) {$\sp_2$};
\draw (.25,-.45) circle [x radius=6mm, y radius=3mm, rotate=45];
 \foreach \x in {1,2,3} {
     \pgfmathparse{int(5-\x)}
     \draw (0,-.6) ++(-\x*35:.5ex) --++(-\x*35-80:1)[fill] circle (.3ex)
         node (u\pgfmathresult)[below] {$\scriptstyle
         u_{\pgfmathparse{int(8-\x)}\pgfmathresult}$}; 
}
\draw[red, rotate=-50] (.9,-.8) arc (-25:-155:.55);
\fill[red] (0.05,-.7) circle (.5ex);
\fill[red] (.85,-1.2) circle (.5ex);
\node (g3) at (.4, -1.2) {$\sp_3$};
\draw (.95,-.9) circle [x radius=6mm, y radius=3mm, rotate=75];
\fill[red] (1.9,-1.2) circle (.5ex);
\node (g4) at (2.4, -1) {$\sp_4$};
\draw (1.8,-.85) circle [x radius=6mm, y radius=3mm, rotate=-75];
\node at (1.9, -1.57) {$\vdots$};
\node[rotate=45] at (2.2, -1.5) {$\vdots$};
\node[rotate=-45] at (1.6, -1.5) {$\vdots$};
\fill[red] (3,-.5) circle (.5ex) node[anchor=base] (g4) {};
\node (g5) at (3.5, 0) {$\sp_5$};
\draw (2.6,-.35) circle [x radius=7mm, y radius=3mm, rotate=-25];
\node at (3, -.87) {$\vdots$};
\node[rotate=45] at (3.3, -.7) {$\vdots$};
\node[rotate=90] at (3.4, -.4) {$\vdots$};
\node (g6) at (-1.5, -.2) {$\sp_6$};
\fill[red] (-1.5,-.6) circle (.5ex) node[anchor=base] (g6) {};
\draw (-1.27,-.6) circle [x radius=5mm, y radius=2.5mm];
\node at (-2, -.6) {$\cdots$};
\node[rotate=-45] at (-1.9, -.3) {$\cdots$};
\node[rotate=45] at (-1.9, -.9) {$\cdots$};
\node[rotate=30] at (-1.3, -1.5) {$\cdots$};
\node[rotate=70] at (-.5, -2.2) {$\cdots$};
 \end{tikzpicture}}}\,.
\end{align*}
Furthermore this construction is reversible: $\sp_1$ is the root of $f$ without
the vertices, using the order on $U$ one obtains a
deterministic way of ordering each underlying $\SP$-structures from vertices.
Finally by induction on the height of $f$ we show that the resulting sequence of
$\SP$-structures satisfies (\ref{eq:pf_condition}).
%
\end{preuve}
\ \\[-10pt]
\begin{minipage}[c]{.66\linewidth}
    \color{fontcolor}
\begin{ex}
    \label{ex:bij_pf_rooted_trees_set_trees}
  Let $(23\mid \cdot \mid 5 \mid 16 \mid \cdot \mid 4\mid \cdot)$ be a (usual)
  parking function on $[6]$. 
  We associate the order $2<_q 3 <_q 5 <_q 1 <_q 6 <_q 4$ from the parking
  function (using the natural order on $[6]$). The bijection described in
  the proof of the Theorem is a generalization of the
  \textsc{Foata-Riordan} bijection \cite{foata1974mappings}.
  The forest of $\F[6]$ associated to this construction is:
\end{ex}
\end{minipage} \hfill
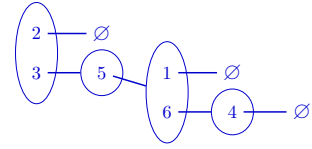
\begin{minipage}[c]{.36\linewidth}
\begin{align*}
\scalebox{.7}{\newcommand{\nna}{\node (a) {$2$};}%
\newcommand{\nnb}{\node (b) {$3$};}%
\newcommand{\nnc}{\node (c) {$5$};}%
\newcommand{\nnd}{\node (d) {$1$};}%
\newcommand{\nne}{\node (e) {$6$};}%
\newcommand{\nnea}{\node (ea) {$\emptyset$};}%
\newcommand{\nneb}{\node (eb) {$\emptyset$};}%
\newcommand{\nnec}{\node (ec) {$4$};}%
\newcommand{\nned}{\node (ed) {$\emptyset$};}%
\begin{tikzpicture}[auto]
\matrix[ampersand replacement=\&, column sep=.75cm, row sep=.2cm]{
  \nna \& \nnea \\
  \nnb \& \nnc \& \nnd \& \nneb \\
       \&      \& \nne \& \nnec \& \nned\\
};
\node[draw,inner sep=2pt,ellipse,fit=(a) (b)] {};
\node[draw,inner sep=2pt,ellipse,fit=(c)](gc) {};
\node[draw,inner sep=2pt,ellipse,fit=(d) (e)] (gde) {};
\node[draw,inner sep=2pt,ellipse,fit=(ec)] (gec) {};
\path[thick] (a) edge (ea)
             (b) edge (gc)
             (c) edge (gde)
             (d) edge (eb) 
             (e) edge (gec)
             (ec) edge (ed);
\end{tikzpicture}}
\end{align*}
\end{minipage}
\\[10pt]

From the species theory this theorem immediately yields the following equalities
of exponential generating series computable by \emph{Lagrange inversion}:
\begin{coro}
    \label{cor:enum_pf_like}
    Let $\SP^\blacktriangleright(t)$ and $\T_\SP(t)$ be respectively the
    exponential generating series of $\SP^\blacktriangleright$ and $\T_\SP$,
    \begin{align*}
    \SP^\blacktriangleright(t) = \T_\SP(t) \qquad\text{and}\qquad
    \T_\SP(t) = \SP(t \T_\SP(t))\,.
    \end{align*}
\end{coro}
\section{Applications on some species}
    \subsection{Linear order species and chained hash tables with linked lists}
We recall the exponential generating series of the linear order:
    $\linSp(t) = (1-t)^{-1}$.
So thanks to the Corollary \ref{cor:enum_pf_like}, the exponential generating
series of $\T_\linSp \simeq \linSp^\btr$ is:
\begin{align*}
    \T_\linSp(t) = \frac{1 - \sqrt{1 - 4t}}{2t}\,. \tag*{\oeis{A001761}}
\end{align*}
\indent
In Example \ref{ex:lin_order_pf} we listed all possible
\emph{hash tables with linked lists} (using parking functions as hash map) on a
set of key $U = [0], [1], [2]$ and $[3]$.
This formula suggests an isomorphism with labeled binary trees; considering
a parking function as a \emph{staircase walk} where
\textit{tread}/horizontal steps are decorated by linear orders of same length
(see \cite[\S2]{priezvirmaux}), the classical bijection between Dyck
paths and binary trees gives a bijection between
$\linSp^\blacktriangleright$-structures and labeled binary trees (when we
create a node from a tread, we associate \textit{via} the bijection the label
to decorate it).

The previous construction (of tree-like structure in the
proof of Theorem \ref{thm:iso}) completed by the bijection between linear orders
and decreasing binary trees provides a bijection between
$\linSp^\blacktriangleright$-structures and ``\textit{labeled incomplete
ternary trees on n vertices in which each left and middle child have a larger
label than their parent}'' (see \oeis{A001761}).
\subsection{Partition species}
The \emph{species of partitions} is defined by the equation:
%
    $\partSp := \setSp(\setSp_+)$ with $\setSp_+$ the
    restriction of $\setSp$ to non-empty sets.
%
Thanks to the corollary \ref{cor:enum_pf_like}, the generating series of
$\partSp^\blacktriangleright$ satisfies the equation:
\begin{align*}
    \partSp^\blacktriangleright(t) = \exp(e^{t\, \partSp^\blacktriangleright(t)}
    - 1) = 1 + t + 4\frac{t^2}{2!} + 29\frac{t^3}{3!} + 311\frac{t^4}{4!} +
    4447\frac{t^5}{5!} + \cdots \tag*{\oeis{A030019}}
\end{align*}
and the construction (proof of the Theorem \ref{thm:iso}) yields the bijection
between $\partSp^\blacktriangleright$-structures and `` spanning trees in the
complete hypergraph on $n$ vertices'' (see \oeis{A030019}).
\subsection{Composition species}
The \emph{species of compositions} (or \textit{ballots}) is defined by
the equation:
%
    $\compSp := \oneSp + \setSp_+ \cdot \compSp$ with $\oneSp$ the
    neutral species ($\oneSp[U] = \{U\}$ \textit{iff} $U = \emptyset$ and it is
    empty in otherwise).
The generating series of the $\compSp^\blacktriangleright$-structures satisfies:
\begin{align*}
    \compSp^\blacktriangleright(t) = \frac{1}{2 - \exp(t\,
    \compSp^\blacktriangleright(t))} = 1 + t + 5\frac{t^2}{2!} +
    46\frac{t^3}{3!} + 631\frac{t^4}{4!} + 11586\frac{t^5}{5!} + \cdots
    \tag*{\oeis{A052894}}
\end{align*}
\subsection{Subsets species}
The \emph{subsets species} is $\mathpzc{S} := \setSp \cdot \setSp$. The
construction of the Theorem \ref{thm:iso} yields a bijection with
\emph{ditrees} (see \oeis{A097629}) and the generating series satisfies $s(t) =
\exp(2t\, s(t))$.
\subsection{$k$-ary trees}
Let $\F_k$ be the species of the $k$-ary trees ($\F_k := \oneSp + \singSp
\cdot (\F_k)^k$ with $k \seq 1$). The construction defines a bijection between
$\F_k^\btr$-structures and $\F_{k+1}$-structures (\textit{Note}:
$\F_1^\btr = \linSp^\btr$). From \cite{beineke1971number}, there also exists an
isomorphism with labeled dissections of a $k$-ball.

\section{Perspective}
    The bijection can be extended between \emph{parking-like structures} over $\SP$
associated to the any non-decrea\-sing linear function $\chi : m \mapsto a\cdot
m + b$ and \emph{tree-like structures} associated to the species $\T_\SP^{a,b}$
defined as $\SP(\singSp \cdot \T_\SP^{a,0})^{a+b}$ (with the root have
$a+b$ ordered edges and others nodes have $a$ ordered edges). This extension is
a generalization of the bijective proof \cite[Theorem 1]{yan2001generalized}.

%
%
%
\small
\bibliographystyle{alpha}
\bibliography{biblio}

\end{document}